\documentclass[10pt,twoside]{article}
\usepackage{latexsym}
\usepackage{mathrsfs}
\usepackage{dsfont}

\usepackage{graphicx}
\usepackage{arydshln}
\def\thebibliography#1{\center{\bf\normalsize References}\list
 {[\arabic{enumi}]}{\settowidth\labelwidth{[#1]}\leftmargin\labelwidth
 \advance\leftmargin\labelsep
 \usecounter{enumi}}
 \def\newblock{\hskip .11em plus .33em minus .07em}
 \sloppy\clubpenalty4000\widowpenalty4000
 \sfcode`\.=1000\relax}

\def\proof{ {\bf Proof.\hspace{0.19cm}}}
\def\endproof{\hfill$\Box$\vspace{0.4cm}}

\def\disp{\displaystyle}
\newcommand{\eqref}[1]{$(\ref{#1})$}
\newcommand{\refeq}[1]{$(\ref{#1})$}
\newcommand{\thb}[1]{{\rm (#1)}}
\def\ol{\overline}
\def\Ga{\alpha}
\def\Gb{\beta}

\def\Gl{\lambda}

\def\GP{\Phi}
\def\mcL{{\mathscr L}}
\def\mcM{{\mathscr M}}

\def\mcU{{\mathscr U}}

\def\bX{\ol{X}}

\def\bb{\bar b}

\def\bu{\bar u}

\def\tiA{\widetilde A}

\def\tiC{\tilde C}

\def\tiX{\widetilde  X}
\def\tiY{\widetilde  Y}

\def\tiu{\tilde u}

\def\tix{\tilde x}

\def\tiGP{\widetilde \GP}
\def\tiPsi{\widetilde \Psi}

\def\mm{\mbox{ }}

\def\bb{\hspace{18pt}}

\def\qq{\qquad}
\def\q{\quad}

\def\eqon{ \mm {\rm on } \mm }

\def\all{  \mm \forall \mm }

\newcommand{\doo}[2]{{d {#1} \over d {#2}}}

\newcommand{\ddo}[3]{{d {^{#1}#2} \over d {#3^{#1}}}}

\def\rank{\mm{\rm rank}\mm}

\def\IC{\mathds{C}}
\def\IR{{\rm I\hspace{-0.90mm}R}\mbox{}}

\def\defeq{\stackrel{\triangle}{=}}

\begin{document}

\title{
{\bf  Necessary and Sufficient Conditions for Distinguishability of
Linear Control Systems}
\thanks{This work was supported in part by NSFC (No.
61074047 and 10831007), and 973 Program (No. 2011CB808002).} }
\author{Hongwei Lou\thanks{School of Mathematical Sciences, and LMNS, Fudan
University, Shanghai 200433, China (Email:
\texttt{hwlou@fudan.edu.cn}). } }

\date{}
\maketitle
\begin{quote}
\footnotesize {\bf Abstract.} Distinguishability takes a crucial
rule in studying observability of hybrid system such as switched
system. Recently, for two linear systems, Lou and Si gave a
condition not only necessary but also sufficient to the
distinguishability of linear systems. However, the condition is not
easy enough to verify. This paper will give a new equivalent
condition which is relatively easy to verify.


\end{quote}
\normalsize

\newtheorem{Definition}{\bb Definition}[section]
\newtheorem{Theorem}[Definition]{\bb Theorem}
\newtheorem{Lemma}[Definition]{\bb Lemma}
\newtheorem{Corollary}[Definition]{\bb Corollary}
\newtheorem{Proposition}[Definition]{\bb Proposition}
\newtheorem{Remark}{\bb Remark}[section]
\newtheorem{Example}{\bb Example}[section]
\newfont{\Bbb}{msbm10 scaled\magstephalf}
\newfont{\frak}{eufm10 scaled\magstephalf}
\newfont{\sfr}{msbm7 scaled\magstephalf}

\def\theequation{1.\arabic{equation}}
\setcounter{equation}{0} \setcounter{Definition}{0}
\setcounter{Remark}{0}
\section{Introduction}
Consider a switched system composed by two time-invariant subsystems
$\ (i=1,2)$:
\begin{equation}\label{E101}
S_i:\q \left\{\begin{array}{l}\disp \doo x t= A_ix(t)+B_iu(t), \\
\disp y(t)= C_ix(t)+G_iu(t),
\end{array}\right.
\end{equation}
where $x(t)\in\IR^n,\ u(t)\in\IR^m$ and $y(t)\in\IR^k$. Naturally,
\begin{equation}\label{E102}
A_i\in \IR^{n\times n}, \q B_i\in \IR^{n\times m}, \q C_i\in
\IR^{k\times n}, \q G_i\in \IR^{k\times m}.
\end{equation}

Switched system is an important case of hybrid systems. When we
consider the observability of switched system composed by
time-invariant subsystems such as system \eqref{E101},
distinguishability takes a crucial rule (see \cite{C-S},
\cite{L-S}). Among the references about
observability/distinguishability of hybrid system, we would like to
refer the readers to the papers \cite{B-B-B-S}, \cite{B-E},
\cite{B-E2}, \cite{B-F-M}, \cite{B-P}, \cite{F-G}, \cite{O-H-T},
\cite{O-W}, \cite{S-B-P}, \cite{V-C-S} and \cite{V-C-S-S}.

In \cite{V-C-S-S}, the authors got a necessary and sufficient
condition for distinguishability of two linear automation  systems
(i.e. $B_1=B_2=0$, $G_1=G_2=0$). However, as pointed out by the
authors of \cite{V-C-S-S}, for non-automation system, the input
plays a crucial role and the distinguishability of two linear
systems becomes very difficult. Recently, in \cite{L-S}, the authors
gave a definition of distinguishability for linear non-automation
systems(see Definition \ref{T101} below), and yielded a necessary
and sufficient condition for distinguishability of two linear
systems.

\begin{Definition} \thb{\textbf{distinguishability}}\label{T101}
Systems $S_1$~and~$S_2$ are said to be \textbf{distinguishable}  on
$[0,T]$,  if for any non-zero
$$
(x_{10},x_{20},u(\cdot))\in \IR^n\times \IR^n \times L^1(0,T;\IR^m),
$$
the corresponding outputs $y_1(\cdot)$ and $y_2(\cdot)$ can not be
identical to each other on $[0,T]$.
\end{Definition}

To study the distinguishability of two systems, some auxiliary
concepts of distinguishability was also introduced in \cite{L-S}:
\begin{Definition}\label{T202} Given $T>0$. Let $\mcU \subseteq L^1(0,T;R^m)$ be a
function space. We say that $S_1$ and $S_2$ are \textbf{$\mcU$ input
distinguishable} on $[0,T]$ if for any non-zero
$$
(x_{10},x_{20},u(\cdot))\in R^n\times R^n \times \mcU,
$$
the outputs $y_1(\cdot)$ and $y_2(\cdot)$ can not be identical to
each other on $[0,T]$.

Especially, when $\mcU$ is the set of polynomial function class, the
set of analytic function class and the set of smooth function  class
$C^\infty([0,T];\IR^m)$, then the corresponding distinguishability
is called  ``polynomial input distinguishability",  ``analytic input
distinguishability" and ``smooth input distinguishability", etc.
\end{Definition}

Denote
\begin{equation}\label{F1}
A=\pmatrix{A_1 & 0\cr 0 & A_2},\, B=\pmatrix{B_1 \cr B_2},\,
C=\pmatrix{C_1 & -C_2},\, G=G_1-G_2
\end{equation}
and
\begin{equation}\label{F2}
X_0=\pmatrix{x_{10}\cr x_{20}}, \,  Y(\cdot)=y_1(\cdot)-y_2(\cdot).
\end{equation}
Then the distinguishability of  $S_1$ and $S_2$ on $[0,T]$ is
equivalent to that for the following system:
\begin{equation}\label{E104}
S:\q \left\{\begin{array}{l}\disp \doo X t= A X(t)+Bu(t), \\
\disp X(0)=X_0;\\
\disp Y(t)= CX(t)+Gu(t),
\end{array}\right.
\end{equation}
$(X_0,u(\cdot))\ne 0$ implies $Y(\cdot)\not\equiv 0$ on $[0,T]$.

It was proved in \cite{L-S} that
\begin{Theorem}\label{T103} The
distinguishability of  $S_1$ and $S_2$ on $[0,T]$ is  equivalent to
that $S_1$ and $S_2$ are analytic input distinguishable.  Moreover,
it is equivalent to that the infinite dimensional equation
{\footnotesize
\begin{equation}\label{E103} \hspace{-8mm}\mcM \Gb\equiv  \pmatrix{C
& G & 0 & 0 &  0 & \ldots \cr
               CA & CB & G & 0 & 0 & \cdots\cr
               CA^2 & CAB & CB & G & 0 & \cdots\cr
               CA^3 & CA^2B & CAB & CB & G& \cdots \cr
\vdots&\vdots&\vdots&\vdots&\vdots&\ddots}\pmatrix{\Gb_1\cr \Gb_2\cr
\Gb_3 \cr \Gb_4\cr  \Gb_5\cr \vdots }=0. \end{equation} } admits
only trivial solution. Consequently, it is independent of $T$.
\end{Theorem}

The disadvantage of Theorem \ref{T103} is that whether   equation
\refeq{E103} admitting only trivial solution is not easy to verify.
In this paper, we will seek for an equivalent condition which can be
verified easier.

\def\theequation{2.\arabic{equation}}
\setcounter{equation}{0} \setcounter{Definition}{0}
\setcounter{Example}{0}
\section{Properties of Differential Operator $D$ and Laplace Transform
$\mcL$}

We recall the notions of differential operator $D$ and Laplace
transform $\mcL$ and list some useful properties of them. For the
cause of notation simplicity, in this section, matrices $A, B,C,G$
and integers $n,m,k$ and etc., can be different from that in other
sections.

Denote by $\IC$ the space of complex numbers. Let
$$
P(\Gl)=a_n\Gl^n+a_{n-1}\Gl^{n-1}+\ldots+a_1\Gl+a_0
$$
be a polynomial, where $a_k\in \IC$ and $\Gl\in \IC$. For smooth
vector valued function $f:[0,T]\to \IC^k$, define
\begin{equation}\label{E201}
P(D)f (t)=a_n\ddo n {f(t)} t  +a_{n-1}\ddo {n-1}{f(t)} t
+\ldots+a_1f^\prime (t)+a_0f(t).
\end{equation}

It is well known that
\begin{Lemma}\label{T201} Let
$f,g:[0,T]\to \IC^n$ be two  smooth vector valued functions on
$[0,T]$, $P(\cdot),Q(\cdot)$ are two polynomials, $\Ga,\Gb\in \IC$
and $\Gl\in \IC$ are two complex constants. Then
\begin{equation}\label{E202}
\Big(\Ga P(D)+\Gb Q(D)\Big)f (t)=\Ga P(D)f(t)+\Gb Q(D)f (t),
\end{equation}
\begin{equation}\label{E203}
P(D)\Big(\Ga f (t)+\Gb g(t)\Big)=\Ga P(D)f(t)+\Gb P(D)g (t),
\end{equation}
\begin{equation}\label{E204}
\Big(P(D)Q(D)\Big)f
(t)=P(D)\Big(Q(D)f(t)\Big)=Q(D)\Big(P(D)f(t)\Big),
\end{equation}
\begin{equation}\label{E205} P(D)\Big(e^{\Gl t} f (t)\Big)=e^{\Gl t}
P(D+\Gl)f(t).
\end{equation}
\end{Lemma}

For a function $f(\cdot)\in L^1_{loc}[0,+\infty)$, the Laplace
transform of $f(\cdot)$ is defined by
$$
F(s)\equiv \mcL\Big(f(\cdot)\Big)(s)\defeq
\int^{+\infty}_0e^{-st}f(t)\, dt, \qq s>0.
$$

It is well known that Laplace transform can be defined for many
functions and even for generalized functions such as $\delta$
function. If there exist $M_1,M_2>0$ such that
$$
|f(t)|\leq M_1e^{M_2t}, \qq\all t>0,
$$
then $\mcL\Big(f(\cdot)\Big)(s)$ is well defined for $s\in
(M_2,+\infty)$. Moreover, $f(\cdot)$ has the form
$$
f(t)=e^{\Gl_1 t}P_1(t)+e^{\Gl_2 t}P_2(t)+\ldots+e^{\Gl_n t}P_n(t)
$$
with $\Gl_k\in \IC$ and $P_k(\cdot)$ being polynomial
($k=1,2,\ldots$) if and only if $\mcL(f)$ is a proper rational
function.

\def\theequation{3.\arabic{equation}}
\setcounter{equation}{0} \setcounter{Definition}{0}
\setcounter{Example}{0}
\section{Main Results} Now we consider the necessary and sufficient conditions for
distinguishability. Let $A,B,C,G$ be defined as in \S 1. By the
discussions of \cite{L-S}, we know that if $S_1$ and $S_2$ are not
distinguishable, then they are not analytic input  distinguishable.
More precisely, there exists a pair $\disp (X_0,u(\cdot)) $ such
that
\begin{equation}\label{E301}
(X_0,u(\cdot))\ne 0,
\end{equation}
\begin{equation}\label{E303}
Ce^{At}x_0+C\int^t_0e^{A(t-s)}Bu(s)\, ds+Gu(t)=0, \qq\all t\geq 0,
\end{equation}
and
\begin{equation}\label{E304}
u(t)=\sum^\infty_{j=0}{\Ga_j\over j!} t^j, \qq t\in [0,+\infty),
\end{equation}
with
\begin{equation}\label{E305}
|\Ga_j|\leq M^{j+1}, \qq\all j=0,1,\ldots
\end{equation}
for some $M>0$.

One can see that if $u(\cdot)$ satisfies
\refeq{E304}---\refeq{E305}, then
$$
|u(t)|\leq Me^{Mt}, \qq\all t\geq 0
$$
and therefore $\mcL\Big(u(\cdot)\Big)(s)$ can be defined for any
$s>M$.

A crucial property we will prove in the following is that
\begin{Lemma}\label{T301}
If $S_1$ and $S_2$ are not distinguishable, then  we can find a pair
$(\bX_0,\bu(\cdot))$ satisfying \refeq{E301}---\refeq{E303} with
\begin{equation}\label{E305A}
\bu(\cdot)=e^{\Gl_1 t}P_1(t)+e^{\Gl_2 t}P_2(t)+\ldots+e^{\Gl_q
t}P_q(t),
\end{equation}
where $\Gl_i\in \IC$ and $P_i(\cdot)$ are vector-valued polynomials
\thb{$i=1,2,\ldots,q$}.
\end{Lemma}
\proof Since $S_1$ and $S_2$ are not distinguishable,  by the
results of \cite{L-S}, there exists a pair $(X_0,u(\cdot))$
satisfies \refeq{E301}---\refeq{E305}. Denote
$$
\GP(s)=\mcL\Big(Ce^{A\cdot}\Big)(s), \q \Psi(s)=\GP(s)B, \q
U(s)=\mcL\Big(u(\cdot)\Big)(s).
$$
Then every element in matrices $\GP$ and $\Psi+G$ are rational
functions.

Consider the Laplace transform of \refeq{E303}, we have
\begin{equation}\label{E306}
\GP(s)X_0+\Big(\Psi(s)+G\Big)U(s)=0.
\end{equation}
Let
$$
r=\rank \Big(\Psi(s_0)+G\Big)=\max_{s\in [0,+\infty)} \rank
\Big(\Psi(s)+G\Big).
$$

\textbf{Case 1:} $r=m$. There exist $1\leq j_1< j_2<\ldots<j_r\leq
k$ such that the matrix $\tiPsi(s)$ composed by $j_i$-th rows
($i=1,2,\ldots r$) of $\Psi(s)+G$ is invertible at $s=s_0$. Then
since every element in  $\tiPsi(\cdot)$ are rational functions, the
determinant of $\tiPsi(\cdot)$ is a rational function and is not
identical to zero. Consequently, $\tiPsi(s)$ is invertible on
$[0,+\infty)$ except for finite points. Let $\tiGP(s)$ be the matrix
composed by $j_i$-th ($i=1,2,\ldots r$) rows of $\GP(s)$. Then
$$
U(s)=\tiPsi(s)^{-1}\GP(s)X_0, \qq s\in [0,+\infty).
$$
Thus, every element in $U(s)$  are rational functions. Moreover,
noting that the inverse Laplace transform of a (non-zero) polynomial
is the linear combination of $\delta$ function and its derivatives,
elements in $U(s)$ must be proper rational functions since
$u(\cdot)$ is analytic. Thus, $u(\cdot)$ has the form \refeq{E304}.
Therefore, in this case, we can get our result by choosing
$\bX_0=X_0$ and $\bu(\cdot)=u(\cdot)$.

\textbf{Case 2:} $r<m$. Let $\bX_0=0$. We will prove that there
exists a $\bu(\cdot)\ne 0$ such that $(\bX_0,\bu(\cdot))$ satisfies
\refeq{E303} and \refeq{E305A}.

Let $1\leq j_1< j_2<\ldots<j_r\leq k$ satisfy that the matrix
$\tiPsi(s)$ composed by $j_i$-th rows ($i=1,2,\ldots r$)  of
$\Psi(s)+G$ has full row rank when $s=s_0$. Then $\tiPsi(s)$ has
full row rank except for finite points. Moreover, the equation
\begin{equation}\label{E307}
\Big(\Psi(s)+G\Big)V(s)=0
\end{equation}
is equivalent to
\begin{equation}\label{E308}
\tiPsi(s)V(s)=0.
\end{equation}

Without loss of generality, suppose that
$$
\tiPsi(s)=\Big( \tiPsi_1(s) \q \tiPsi_2(s)\Big),
$$
where $\tiPsi_1(s)$ is an $r\times r$ matrix-valued function such
that $\tiPsi_1(s)$ is invertible at $s=s_0$. Then since elements in
$\tiPsi_1(s)$ are rational functions, $\tiPsi_1(s)$ is invertible
except for finite points. Consequently, it is easy to see that
\refeq{E308} admits a solution $V_1(\cdot)$ with
$$
V_1(s)=\pmatrix{Q_1(s)\cr  \vdots \cr Q_r(s)\cr 1\cr \vdots\cr 1},
$$
where $Q_j(s)$ ($j=1,\ldots r$) are rational functions. Choosing $J$
large enough and letting
$$
V(s)={V_1(s)\over s^J},
$$
we get a nontrivial solution $V(\cdot)$ of \refeq{E308}(or
\refeq{E307}, equivalently) such that every element of $V(\cdot)$
are proper rational functions. Consequently,
$$
\bu(\cdot)=\mcL^{-1}(V(\cdot))\ne 0
$$
is well-defined and $\bu(\cdot)$ has the form \refeq{E305A}.
Moreover,
$$
C\int^t_0e^{A(t-s)}B\bu(s)\, ds+G\bu(t)=0, \qq\all t\geq 0.
$$
That is $(\bX_0,\bu(\cdot))$ satisfies \refeq{E303}.

We get the proof.
\endproof

Using the properties of the differential operator $D$, we can go
further.

\begin{Lemma}\label{T302}
If $S_1$ and $S_2$ are not distinguishable, then we can find a pair
$(\tiX_0,\tiu(\cdot))\ne 0$ satisfying \refeq{E303} and
\begin{equation}\label{E309}
\tiu(\cdot)=e^{\Gl  t}\xi,
\end{equation}
where $\Gl \in \IC$ and $\xi\in\IC^m$.
\end{Lemma}
\proof By Lemma \ref{T301}, there exists a pair $(\bX_0,\bu(\cdot))$
satisfying \refeq{E301}---\refeq{E303} and \refeq{E305A}.

\textbf{Case 1:}  $\bu(\cdot)\equiv 0$. Then let
$\tiu(\cdot)=\bu(\cdot)$, we get the conclusion.

\textbf{Case 2:}  $\bu(\cdot)\not\equiv 0$. Then
$$
\bu(t)=e^{\Gl_1 t}P_1(t)+e^{\Gl_2 t}P_2(t)+\ldots+e^{\Gl_q t}P_q(t),
$$
where
$$
P_i(t)=\xi_{p_i,i}t^{p_i}+\xi_{p_i-1,i}t^{p_i-1}+\ldots+\xi_{1,i}t+\xi_{0,i},\q
i=1,\ldots, q.
$$
$$
p_i\geq 0, \q \xi_{p_i,i}\ne 0,
$$
and
$$
\Gl_i\ne \Gl_j, \qq i\ne j.
$$
Let $(X(\cdot),Y(\cdot))$ be the solution of \refeq{E104}
corresponding to $(\bX_0,\bu(\cdot))$. Then
$$
Y(t)\equiv  CX(t)+G\bu(t)\equiv 0
$$
since $(\bX_0,\bu(\cdot))$ satisfies \refeq{E303}.

Let
$$
Q(\Gl)=(\Gl-\Gl_1)^{p_1}(\Gl-\Gl_2)^{p_2+1}\ldots
(\Gl-\Gl_q)^{p_q+1},
$$
$$
\tiu(t)=Q(D)\bu(t)
$$
and
$$
\tiX(t)=Q(D)X(t), \q \tiY(t)=Q(D)Y(t).
$$
We have $\tiY\equiv 0$,
$$
\doo \tiX t= A\tiX(t)+B\tiu(t),
$$
and
$$
\tiY(t)= C\tiX(t)+G\tiu(t).
$$
That is, $(\tiX(\cdot),\tiY(\cdot))$ is the solution of \refeq{E104}
corresponding to $(\tiX_0,\tiu(\cdot))$ for some $\tiX_0\in
\IR^{2n}$. In other words, $(\tiX_0,\tiu(\cdot)) $  satisfies
\refeq{E303}.

Finally, \refeq{E301} follows from
\begin{eqnarray*}\disp
\tiu(t) &=& Q(D)\bu(t) \\
\disp &=& \sum^q_{i=1} e^{\Gl_i t}Q(D+\Gl_i) P_i(t) \\
\disp &=& e^{\Gl_1 t}Q(D+\Gl_1) P_1(t) \\
\disp &=& e^{\Gl_1 t}(D+\Gl_1-\Gl_2)^{p_2+1}\ldots (D+\Gl_1-\Gl_q)^{p_q+1} D^{p_1}P_1(t) \\
\disp &=& (p_1)!(\Gl_1-\Gl_2)^{p_2+1}\ldots
(\Gl_1-\Gl_q)^{p_q+1}e^{\Gl_1 t}\xi_{p_1,1}\\
\disp &\ne & 0.
\end{eqnarray*}
\endproof

\begin{Corollary}\label{T304}
If $S_1$ and $S_2$ are $0$-th  polynomial input distinguishable,
then they are  polynomial input distinguishable.
\end{Corollary}
\proof Subsystems $S_1$ and $S_2$ are $0$-th  polynomial input
distinguishable means that for any $(x_{10},x_{20},u(\cdot))\in
R^n\times R^n \times \mcU$ with $u(\cdot)\equiv \xi\in \IR^m$, the
outputs $y_1(\cdot)$ and $y_2(\cdot)$ can not be identical to each
other on $[0,T]$.

If  $S_1$ and $S_2$ are not polynomial input distinguishable, then
there exists $(X_0,u(\cdot))$ $\ne 0$ such that \refeq{E303} holds
with $u(\cdot)$ being a polynomial. Then using the method we
constructed $\tiu(\cdot)$ from $\bu(\cdot)$ in the proof of Lemma
\ref{T302}, we can construct a pair $(\tiX_0,\tiu(\cdot))\ne 0$
satisfying \refeq{E303} with $\tiu(\cdot)$ being a constant vector.
This means $S_1$ and $S_2$ are not $0$-th  polynomial input
distinguishable.
\endproof

By Corollary \ref{T304}, the necessary and sufficient condition for
$0$-th  polynomial input distinguishable and that for
 $k$-th  polynomial input
distinguishable are equivalent. Thus, by Theorem 3.1 of \cite{L-S},
we can see that for any $p\geq 0$, the matrix {\footnotesize $$
\hspace{-8mm} \pmatrix{C & G & 0 & \ldots & 0 \cr
               CA & CB & G & \cdots &  0 \cr
\vdots&\vdots&\vdots&\ddots&\vdots\cr
               CA^{p+2} & CA^{p+1}B & CA^pB & \cdots & G\cr
\vdots&\vdots&\vdots&\vdots&\vdots}$$ } has full column rank if and
only if
$$ \hspace{-8mm}\pmatrix{C & G \cr
               CA & CB \cr
               CA^2 & CAB \cr
               CA^3 & CA^2B \cr
\vdots&\vdots}
$$
has full column rank. Thus, it follows from Cayley-Hamilton's
theorem, they are both equivalent to that
$$ \hspace{-8mm}\pmatrix{C & G \cr
               CA & CB \cr
               CA^2 & CAB \cr
\vdots&\vdots\cr
               CA^{2n} & CA^{2n-1}B }
$$
has full column rank.

Now, we state our main result.
\begin{Theorem}\label{T303}
Systems $S_1$ and $S_2$ are distinguishable if and if only for any
$\Gl\in \IC$,
\begin{equation}\label{E310}
\mcM_\Gl\equiv \pmatrix{C & G \cr
        C(A-\Gl I) & CB\cr
        C(A-\Gl I)^2 & C(A-\Gl I)B\cr
        \vdots & \vdots\cr
        C(A-\Gl I)^{2n} & C(A-\Gl I)^{2n-1}B}
\end{equation}
has full column rank.
\end{Theorem}
\proof (i) Suppose that $S_1$ and $S_2$ are distinguishable. Let
$\Gl\in \IC$. Consider
\begin{equation}\label{E311}
\left\{\begin{array}{l} \disp\doo {\tiX (t)}t=(A-\Gl I)\tiX(t)+B\tiu(t),\\
\disp \tiX(0)=\tiX_0,\\
\disp \tiY(t)=C\tiX(t)+G\tiu(t).
\end{array}\right.
\end{equation}
We claim for any $\tix\in \IC^{2n}$ and $\xi\in \IC^m$,
$(\tiX_0,\xi)\ne 0$, the solution of \refeq{E311} corresponding to
$\tiX_0$ and
$$
\tiu(t)\equiv \xi
$$
satisfies
$$
\tiY(t)\not\equiv 0, \qq \eqon [0,+\infty).
$$
In other words, $(A_1-\Gl I, B_1, C_1,G_1)$ and  $(A_2-\Gl I, B_2,
C_2,G_2)$ are $0$-th polynomial input distinguishable\footnote{Here
variables are complex, but this can be treated similarly to the case
of that only real variables are concerned.}.

If it is not the case, then we have $(\tiX_0,\xi)\ne 0$ such that
the corresponding $\tiY(\cdot)$ equals to zero identically.

Let
$$
X(t)=e^{\Gl t} \tiX(t), \q Y(t)=e^{\Gl t} \tiY(t),
$$
Then $(X(\cdot),Y(\cdot))$ solves \refeq{E104} with
$$
X_0=\tiX_0, \q u(t)=e^{\Gl t} \tiu(t).
$$
Since
$$
Y(t)=e^{\Gl t} \tiY(t)=0,
$$
by considering the real part or imaginary part of
$X_0,u(\cdot),X(\cdot)$ and $Y(\cdot)$, one can  easily see that
$S_1$ and $S_2$ are not distinguishable. This is a contradiction.

Similar to Theorem 3.1 of \cite{L-S}, we can get that the $0$-th
polynomial input distinguishable of $(A_1-\Gl I, B_1, C_1,G_1)$ and
$(A_2-\Gl I, B_2, C_2,G_2)$ (in complex variable sense) implies that
$\mcM$ has full column rank.

(ii) Suppose that $S_1$ and $S_2$ are not distinguishable. Then,
Lemma \ref{T302} shows that there is a pair $(X_0,u(\cdot))\ne 0$
satisfying \refeq{E303} and
\begin{equation}\label{E312}
 u(\cdot)=e^{\Gl t}\xi,
\end{equation}
for some $\Gl\in \IC$. This implies that $(A_1-\Gl I, B_1, C_1,G_1)$
and $(A_2-\Gl I, B_2, C_2,G_2)$ are not $0$-th polynomial input
distinguishable. Consequently, $\mcM_\Gl$ has not full column rank.
\endproof

\def\theequation{4.\arabic{equation}}
\setcounter{equation}{0} \setcounter{Definition}{0}
\setcounter{Example}{0}
\section{Generalization} In \S 1, the  state variables are taken values in
$\IR^n$. In fact, we can consider more general cases. That is,
 for subsystem $S_i$ of  \refeq{E101}, suppose that
$$
A_i\in \IR^{n_i\times n_i}, \q B_i\in \IR^{n_i\times m}, \q C_i\in
\IR^{k\times n_i}, \q G_i\in \IR^{k\times m},  \q i=1,2,
$$
where $n_1,n_2,k,m\geq 1$.

Similar to Definition \ref{T101}, we define
\begin{Definition} \label{T401}
Systems $S_1$~and~$S_2$ are said to be \textbf{distinguishable}  on
$[0,T]$,  if for any non-zero
$$
(x_{10},x_{20},u(\cdot))\in \IR^{n_1}\times \IR^{n_2} \times
L^1(0,T;\IR^m),
$$
the corresponding output $y_1(\cdot;x_{10},u(\cdot))$ of $S_1$
(satisfying the initial condition $x(0)=x_{10}$) and
$y_1(\cdot;x_{20},u(\cdot))$ of $S_2$ (satisfying the initial
condition $x(0)=x_{20}$) are not identical to each other on $[0,T]$.
\end{Definition}
We have
\begin{Theorem} \label{T402}
Subsystems $S_1$ and $S_2$ are distinguishable if and only if for
any $\Gl\in \IC$, the matrix
\begin{equation}\label{E401}
\mcM_\Gl\equiv \pmatrix{C & G \cr
        C(A-\Gl I) & CB\cr
        C(A-\Gl I)^2 & C(A-\Gl I)B\cr
        \vdots & \vdots\cr
        C(A-\Gl I)^{n_1+n_2} & C(A-\Gl I)^{n_1+n_2-1}B}
\end{equation}
has full column rank, where
$$
A=\pmatrix{A_1 & 0\cr 0 & A_2}\in \IR^{(n_1+n_2)\times (n_1+n_2)},\q
B=\pmatrix{B_1 \cr B_2}\in \IR^{(n_1+n_2)\times m},
$$
$$
C=\pmatrix{C_1 & -C_2}\in \IR^{k\times (n_1+n_2)},\q G=G_1-G_2\in
\IR^{k\times m}.
$$
\end{Theorem}

At the end of the paper, we prove the following equivalent result.
\begin{Theorem}\label{T501}
Let $\mcM_\Gl$ be defined by \refeq{E401}. Then  $\mcM_\Gl$ has full
column rank or any $\Gl\in \IC$  if and only if for any $\Gl\in \IC$
$$
\pmatrix{C & G\cr A-\Gl I & B}
$$
has full column rank.
\end{Theorem}
\proof Let $\disp \tiC= \pmatrix{C & G},\q \tiA_\Gl=\pmatrix{A-\Gl I
& B\cr 0 & 0}. $ Then it is well known that for any ~$\Gl\in\IC$,
$$ \pmatrix{ \tiC \cr
\tiC \tiA_\Gl    \cr \vdots \cr \tiC\tiA_\Gl^{n_1+n_2+m-1} }
$$
has full column rank if and only if for any ~$\Gl,\mu \in \IC$,
$\disp \pmatrix{ \tiC  \cr \tiA_\Gl-\mu I } $ has full column rank
(see \cite{S-S} for example).

Thus, using Cayley-Hamilton's Theorem, for any~$\Gl\in\IC$,
 $\disp\mcM_\Gl$ has full column rank if and only if for any~$\Gl,\mu
\in \IC$, $$\disp \pmatrix{ C & G \cr A-\Gl I  & B\cr 0 & -\mu I }
$$ has full column rank. While the later is equivalent to that for
any ~$\Gl\in \IC$, $\disp \pmatrix{ C & G \cr A-\Gl I & B } $ has
full column rank. We get the proof.
\endproof

\footnotesize
\vspace{5mm}\footnotesize \ \\

\end{document}